\documentclass[11pt]{article}
\usepackage{amsmath,amssymb,amscd,graphicx, color}
\usepackage{hyperref}
\usepackage{tensor}

\newcommand{\bbR}{{ \mathbb{R} }}

\newcommand{\R}{\mathbb{R}}
\newcommand{\dx}{{\,{\mathrm{d}}x}}

\newtheorem{lemma}{Lemma}

\newtheorem{theorem}{Theorem}

\newtheorem{remark}{Remark}
\numberwithin{equation}{section}

\begin{document}
\title{Some Liouville-type theorems for the stationary 3D magneto-micropolar fluids}
\renewcommand{\baselinestretch}{1.2}
\renewcommand{\theequation}{\thesection.\arabic{equation}}
\author{Jae-Myoung Kim\thanks{Department of Mathematics Education, Andong National University, Andong, Republic of Korea. Email: \tt{jmkim02@andong.ac.kr}}
~and ~Seungchan Ko\thanks{Department of Mathematics, Inha University, Incheon, Republic of Korea. Email:
\tt{scko@inha.ac.kr}}}

\date{}
\maketitle
\begin{abstract}
In this paper, we prove some Liouville-type theorems for the stationary magneto-micropolar fluids under suitable conditions in three space dimensions. We first prove that the solutions are trivial under the assumption of certain growth conditions for the mean oscillations of the potentials. And then we show similar results assuming that the solutions are contained in $L^p(\R^3)$ with $p\in[2,9/2)$. Finally, we show the same result for lower values of $p\in[1,9/4)$ with the further assumption that the solutions vanish at infinity.  
  \noindent{\bf
keywords:} stationary magneto-micropolar equations, Liouville-type
theorem
\end{abstract}

\section{Introduction}
In the present paper, we consider the stationary magneto-micropolar fluid equations in $\mathbb{R}^3$, which consists of the following partial differential equations:
\begin{equation}
\label{mag-mipolar} \left\{
\begin{aligned}
-\Delta u+( u \cdot \nabla) u + \nabla \Pi& =  \chi \nabla \times w + (b \cdot \nabla)  b,  \\
 -\gamma \Delta w+(u \cdot \nabla) w & =   \nabla (\nabla \cdot w) + \chi \nabla \times u -2 \chi w,  \\
-\nu \Delta b+ (u \cdot \nabla) b & = (b \cdot \nabla) u, \\
\nabla \cdot u  & =  \nabla \cdot b\, = 0,
\end{aligned}
\right.
\end{equation}
where $u=(u_1,u_2,u_3)$, $w=(w_1,w_2,w_3)$, $b=(b_1,b_2,b_3)$ and $\Pi$ denote the fluid velocity,
the angular velocity of the rotation of the fluid particles, the
magnetic fields and pressure respectively. The positive constant
$\gamma$ in (\ref{mag-mipolar}) correspond to the angular viscosity,
$\nu$ is the inverse of the magnetic Reynolds number and $\chi$ is
the micro-rotational viscosity. In this paper, we assume that $\gamma$=$\nu$=$\chi$=1 for simplicity. Equation $\eqref{mag-mipolar}_1$ is similar to the classical Navier--Stokes equations, but here it is coupled with equations $\eqref{mag-mipolar}_2$ for $w$ and $\eqref{mag-mipolar}_3$ for $b$. Equation $\eqref{mag-mipolar}_2$ describes the motion in the macro-volumes as they go through micro-rotational effects, represented by the micro-rotational velocity vector $w$. If the fluids have no microstructure, $w$ vanishes and the system \eqref{mag-mipolar} becomes a magneto-hydrodynamics system. Equation $\eqref{mag-mipolar}_3$ is the Maxwell system for the electric field. This model was first introduced by Ahmadi and Shahinpoor \cite{model}. After that, Rojas-Medar \cite{exist_1} proved the local-in-time existence and uniqueness of strong solutions in a bounded domain based on the spectral Galerkin method. Furthermore, Rojas-Medar and Boldrini \cite{exist_2} established the existence of weak solutions to the model \eqref{mag-mipolar} in a bounded domain, and in particular, the uniqueness was also proved for a two-dimensional domain. The existence of global-in-time strong solutions was addressed by Ortega-Torres and Rojas-Medar \cite{exist_3}. 

After Galdi's work in \cite{G11}, Liouville-type problems for the
stationary fluid equations have been extensively studied and there are
a large number of works on the Liouville type-problems even to these dates (see e.g. \cite{CY2013, C14, CW16, KPR15, SW20} and a review paper \cite{Zhang21}). Here, we shall study some Liouville-type results under the assumptions with regard to the potential functions.
We say that $\Phi \in L^1_{\rm{loc}}(\bbR^3;\bbR^{3\times 3})$ is the
potential functions for the vector fields $u \in L^1_{\rm{loc}}(\bbR^3)$,
if $\mathrm{div}\, \Phi = u$. In \cite{S16}, Seregin obtained Liouville-type theorems for the steady-state Navier-Stokes equations under the assumption that the potential $\Phi \in
{\rm{BMO}}(\bbR^3)$ and $u\in L^6(\R^3)$, and in \cite{S18} the integrability condition for the velocity was dropped. After that, very recently, Chae and Wolf \cite{CW19}
showed Liouville-type theorem for the stationary Navier-Stokes
equations under the assumption
$$
\left( \frac{1}{|B_r|}\int_{B_r} |\Phi - \Phi_{B_r} |^s
\dx \right)^{\frac 1s} \lesssim  r^{ \frac 13 - \frac 1s}
\qquad \forall  1 < r < + \infty
$$
for some $3 < s < + \infty$, and similar results were proved for MHD equations in \cite{CKW211-1}. The first theorem of the present paper is the extension of the result of \cite{CKW211-1}. Here, however, we shall adopt a different approach to control the pressure term by introducing an auxiliary function and utilizing it as a test function. In specific, we first aim in this paper to prove the following Liouville-type result.

\begin{theorem}\label{thm1}
    Let $(u,b,w,\Pi)$ be a smooth solution to the equations \eqref{mag-mipolar}. Assume that there exist potentials $\Phi,\Psi, \Upsilon \in C^\infty(\bbR^3 ; \bbR^{3\times 3})$ such that $\nabla \cdot \Phi = u$, $\nabla \cdot \Psi = b$, $\nabla \cdot \Upsilon = w$ and
\[
\left( \frac 1 {| B_r|} \int_{B_r} \big|\Phi - \Phi_{B_r} \big|^s
\dx \right)^{\frac 1s}
        + \left( \frac 1 {| B_r |} \int_{B_r} \big|\Psi - \Psi_{B_r} \big|^s \dx \right)^{\frac 1s}
\]
    \begin{equation}\label{AS}
+\left( \frac 1 {| B_r |} \int_{B_r} \big|\Upsilon -
\Upsilon_{B_r} \big|^s \dx \right)^{\frac 1s}
        \leq C r ^ {\frac 13 - \frac 1s}, \qquad r > 1
    \end{equation}
    for some $3 < s \leq 6$. Then $u \equiv b \equiv w \equiv 0$.
\end{theorem}
\begin{remark}
    In the case of $w\equiv 0$, Theorem~\ref{thm1} reduces to  \cite[Theorem~1.1]{CKW211-1}.
\end{remark}

Later, Zhang et. al. \cite{ZYQ15} proved that if smooth solutions
of the stationary MHD equations are bounded in $L^{\frac{9}{2}}(\bbR^3)$ and have finite Dirichlet integral, then they are also identically zero. After that, Schulz \cite{Sc19} obtained the
Liouville theorem for this equations provided that the smooth solution $(u, b)$ are contained in $L^p(\R^3)\cap {\rm{BMO}}^{-1}(\bbR^3)$ with $p \in(2, 6]$. Recently, Yuan and Xiao \cite{YX20} proved
that if smooth solution $(u, b) \in L^p(\bbR^3)$ with $2 \leq p\leq
\frac{9}{2}$, then $u =b =0$. In this direction, the second objective of this paper is as follows.
\begin{theorem}\label{thm2}
Let $p \in [2, \frac{9}{2})$. Assume that $( u, b, w, \Pi)$ is a
smooth solution to the equations \eqref{mag-mipolar} with $ u, b, w
\in L^p(\bbR^3)$.
 Then $u \equiv b \equiv w \equiv 0$.
\end{theorem}

Furthermore, parallel to the result of Liu and Liu \cite{LL22}, we shall also prove the following theorem.
\begin{theorem}\label{thm3}
Let $p \in [1, \frac{9}{4})$. Assume that $( u, b, w, \Pi)$ is a
smooth solution to the equations \eqref{mag-mipolar} with $ u, b, w
\in L^p(\bbR^3)$ satisfying $ \lim_{|x|\rightarrow\infty}\ u(x) =
\lim_{|x|\rightarrow\infty}\ b(x)=\lim_{|x|\rightarrow\infty}\ w(x)=0$.
 Then $u \equiv b \equiv w \equiv 0$.
\end{theorem}

\begin{remark}\label{rk}
Even if we consider the model with the variable density, that is, density-dependent models, Theorem \ref{thm1}, \ref{thm2} and \ref{thm3} still hold under suitable additional assumptions (see, for example, \cite{LL22}).
\end{remark}

\begin{remark}\label{rkrk}
In the light of the work of Liu and Liu \cite{LL22}, through a similar
approach, we can also obtain the Liouville-type results in Lorentz spaces (see e.g. \cite{KTW17, CJL21}).
\end{remark}

\section{Preliminaries}

In this section, we introduce some notations and auxiliary results which will be used throughout the paper. We denote the ball with center $x_0$ and radius $R$ by $B_R(x_0)$. If $x_0=0$, we simply write $B_R=B_R(0)$. Throughout the paper, the notation $P\lesssim Q$ implies that there exists some constant $C>0$ such that $P\leq CQ$. Also, $C$ denotes a generic positive constant, which may change at each appearance.

Let us also define a family of cut-off functions. For $0<r<r'$, we let $\xi=\xi_{r,r'} \in C^\infty_c(B_{r'})$ be a radially non-increasing scalar function such that
\begin{equation}\label{cutoff}
\xi_{r,r'}(x)=\begin{cases}
1,\quad x\in  B_r,\\
0,\quad x\in  B^c_{r'},
\end{cases}
\end{equation}
with the properties $|\nabla \xi_{r,r'} | < C_1/ (r'-r)$, and $|\nabla^2 \xi_{r,r'} | < C_2/ (r'-r)^2$ for some constant $C_1,$ $C_2>0$.

Next, for a bounded domain $\Omega\subset\R^3$, we consider the following problem: for given $f\in L^p(\Omega)$ with
\begin{equation}\label{mean_zero}
\int_{\Omega}f(x)\dx=0,
\end{equation}
find a vector-valued function $v\in W^{1,p}_0(\Omega)^d$ satisfying
\begin{equation}\label{aux_prob}
\begin{aligned}
\nabla\cdot v&=f,\\
\|\nabla v\|_p&\leq C\|f\|_p
\end{aligned}
\end{equation}
for some constant $C=C(d,p,\Omega)$. For this matter, we have the following theorem which is quoted from \cite{G11}.
\begin{lemma}\cite[Theorem III.3.1]{G11}\label{aux_thm}
Assume that $\Omega$ satisfies the cone condition. Then for given $f\in L^p(\Omega)$ with $1<p<\infty$ satisfying \eqref{mean_zero}, there exists at least one solution for the problem \eqref{aux_prob}.
\end{lemma}

We will use this result to construct an auxiliary function to control the pressure term in the next section. We will also use the following iteration lemma frequently in our analysis.
\begin{lemma}\cite[Lemma 3.1]{G83}\label{key_2}
Let $f(r)$ be a non-negative bounded function on $[r_0,r_1] \subset \bbR_{\geq0}$. Suppose that there exist non-negative constants $A$, $B$, $D$, $E$ and positive numbers $d <b <a$ and a parameter $\theta\in (0, 1)$ such that for any $r_0\leq s <t \leq r_1$,
\[
f(s) \leq \theta f(t) + \frac{A}{(t- s)^a}  + \frac{B}{(t- s)^b} +
\frac{D}{(t- s)^d} + E.
\]
Then we have
\[
f(s) \leq C(a,b,d,\theta)\left[\frac{A}{(t- s)^a}  + \frac{B}{(t-
s)^b} + \frac{D}{(t- s)^d} + E\right].
\]
\end{lemma}

We shall also use the following lemma for the proof of Theorem \ref{thm1}.
\begin{lemma}\label{LEM3}\cite[Lemma~2.1 and Lemma~2.2]{CKW21}
    Suppose that $R>1$ and $f \in W^{1,2}(B_R;\bbR^3)$. For $0 < \rho <R$, we let $\psi \in C^\infty_c(B_R)$ such that $0 \leq \psi \leq 1$ and $|\nabla \psi | \leq C/(R - \rho)$ for some constant $C>0$. Assume further that there exists the potential $F \in W^{2,2}(B_R;\bbR^{3\times3})$ with $\nabla \cdot F = f$ and the growth condition
    \begin{equation*}
    \left( \frac 1 {| B_r |} \int_{B_r} \big|F - F_{B_r} \big|^s \dx \right)^{\frac 1s} \lesssim r^{\frac 13 - \frac 1s}, \qquad r>1
    \end{equation*}
    for some $3 < s \leq 6$. Then there holds
    \begin{equation}\label{C1}
    \| \psi^2 f \|_{L^2(B_R)}^2 \lesssim R^{\frac {11}6 - \frac 1s} \| \psi \nabla f \|_{L^2(B_R)} + R^{\frac {11}3 - \frac 2s} (R- \rho)^{-2}
    \end{equation}
    and
    \begin{equation}\label{C2}
    \| \psi^3 f \|_{L^3(B_R)}^3 \lesssim R \| \psi \nabla f \|_{L^2(B_R)}^{\frac {18}{s+6}} +  R^{4-\frac 3s} (R - \rho)^{-3} +  R \left( (R-\rho)^{-1} \| \psi^2 f \|_{L^2(B_R)} \right)^{\frac {18}{s+6}}.
    \end{equation}
\end{lemma}

\section{Proof of Theorem~\ref{thm1}}

Let $\varphi_R$ be a cut-off function in $C^\infty_c(\bbR^3)$ given
by $\varphi_R=\xi_{\rho,\tau}$ for $1<\frac{R}{2}<\rho<\frac{2}{3}R<R<\tau<2R$. We begin with some estimates for the terms related to $w$. By using the H\"{o}lder and Young's inequality, we note that
\begin{align*}
-\int_{\R^3}\nabla(\nabla\cdot w)(w\varphi_R^2)\dx
&=\int_{\R^3}(\nabla\cdot w)\nabla\cdot(w\varphi_R^2)\dx\\
&=\int_{\R^3}|\nabla\cdot w|^2\varphi^2_R\dx+2\int_{\R^3}(\nabla\cdot w) w\varphi_R\cdot\nabla\varphi_R\dx\\
&\geq\int_{\R^3}|\nabla\cdot w|^2\varphi^2_R\dx-2\bigg|\int_{\R^3}(\nabla\cdot w) w\varphi_R\cdot\nabla\varphi_R\dx\bigg|\\
&\geq\int_{\R^3}|\nabla\cdot w|^2\varphi^2_R\dx-\varepsilon\int_{\R^3}|\nabla\cdot w|^2\varphi^2_R\dx-C(\varepsilon)\int_{\R^3}|w|^2|\nabla\varphi_R|^2\dx.
\end{align*}
Next, using the vector identity $\nabla\times (u\varphi^2_R)=\varphi^2_R\nabla\times u+\nabla \varphi^2_R\times u$ and Young's inequality, we also note that
\begin{align*}
\int_{\R^3}\nabla\times w\cdot u\varphi^2_R\dx+\int_{\R^3}\nabla\times u\cdot w\varphi^2_R\dx
&=\int_{\R^3}w\cdot\nabla\times(u\varphi^2_R)\dx+\int_{\R^3}\nabla\times u\cdot w\varphi^2_R\dx\\
&=2\int_{\R^3}w\varphi^2_R\cdot\nabla\times u\dx+2\int_{\R^3}w\cdot\varphi_R\nabla\varphi\times u\dx\\
&\leq \frac{2}{3}\int_{\R^3}|\nabla u|^2\varphi^2_R\dx+\frac{3}{2}\int_{\R^3}|w|^2\varphi^2_R\dx\\
&\hspace{4mm}+\varepsilon\int_{\R^3}|w|^2\varphi_R^2\dx+C(\varepsilon)\int_{\R^3}|u|^2|\nabla\varphi_R|^2\dx.
\end{align*}

Before proceeding more, we introduce an auxiliary function that is needed to handle the pressure term. We set $\Omega = B_\tau\setminus B_\rho$ and $f=\nabla\cdot(\varphi^2_Ru)$. Note that $\Omega$ stisfies the cone condition (see, e.g., \cite[Remark III.3.4]{G11}). As we know by Green's Theorem that
\[\int_{B_\tau\setminus B_\rho}\nabla\cdot(\varphi^2_Ru)\dx=\int_{B_{\tau}}\nabla\cdot(\varphi^2_Ru)\dx-\int_{B_{\rho}}\nabla\cdot(\varphi^2_Ru)\dx=\int_{\partial B_\tau}\varphi^2_R u\cdot\nu\dx=0,\]
we can apply Lemma \ref{aux_thm} to show the existence of vector-valued function $W_R\in W^{1,p}_0(B_\tau\setminus B_\rho)$ satisfying
\begin{equation}\label{press}
\nabla\cdot W_R=\nabla\cdot(\varphi^2_Ru)\quad{\rm{in}}\,\,B_\tau\setminus B_\rho,
\end{equation}
with
\begin{equation}\label{aux_est}
\|\nabla W_R\|_{L^p(B_\tau\setminus B_\rho)}
\lesssim\|\nabla\cdot(\varphi^2_Ru)\|_{L^p(B_\tau\setminus B_\rho)}
=\|\nabla(\varphi^2_R)\cdot u\|_{L^p(B_\tau\setminus B_\rho)}
\lesssim
\|\nabla\varphi_R\cdot u\|_{L^p(B_\tau\setminus B_\rho)},
\end{equation}
where we have used the fact that $\nabla\cdot u = 0$ in the last equality.

Now, we multiply the equations $\eqref{mag-mipolar}_1$,
$\eqref{mag-mipolar}_2$ and $\eqref{mag-mipolar}_3$ by
$u\varphi^2_R-W_R$, $b\varphi^2_R$ and $w\varphi^2_R$,
respectively and integrate over $\bbR^3$. Then integration by parts with the use of divergence-free conditions yields that
\begin{align*}
&\int_{\R^3}\left(|\nabla u|^2+|\nabla b|^2+ |\nabla w|^2\right)\varphi_R^2\dx+\int_{\R^3}|\nabla\cdot w|^2\varphi^2_R\dx+\int_{\R^3}|w|^2\varphi^2_R\dx\\
&\lesssim\int_{\R^3}\left(|u|^2+|b|^2+|w|^2\right)|\nabla\varphi_R|^2\dx + \int_{\R^3}\left(|u|^2+|b|^2+|w|^2\right)u\cdot \varphi_R\nabla\varphi_R\dx\\
& \hspace{4mm}+\int_{\R^3}\nabla u\cdot\nabla W_R\dx-\int_{\R^3}(u\cdot\nabla)u\cdot W_R\dx-\int_{\R^3}(b\cdot\nabla) b\cdot W_R\dx\\
\end{align*}
\vspace{-12mm}
\begin{equation}\label{main_ineq}
\hspace{-39mm}+\int_{\R^3}\nabla\times w\cdot W_R\dx+\int_{\R^3}(u\cdot B)(B\cdot\nabla)\varphi^2_R\dx,
\end{equation}
where the pressure term vanishes due to the equality \eqref{press}. We shall estimate the terms on the right-hand side of \eqref{main_ineq}. First, by \eqref{aux_est}, H\"{o}lder's inequality and Young's inequalities we have
\[
\int_{\bbR^3} \nabla u\cdot \nabla W_R \,\mathrm{d}x\leq \|\nabla
 u\|_{L^2(B_\tau)}\|\nabla W_R\|_{L^2(B_\tau)}\leq \frac{1}{4}\|\nabla
 u\|^2_{L^2(B_\tau)}+C(\tau-\rho)^{-2}\|u\|^2_{L^2(B_\tau\setminus
 B_\rho)}.
\]
Next, by \eqref{aux_est} and H\"older's inequality, we obtain

\begin{align*}
\int_{\R^3}(u\cdot\nabla)u\cdot W_R\dx&=\sum_{i,j}\int_{B_{\tau}\setminus B_{\rho}}u_i\partial_iu_j(W_R)_j\dx=-\sum_{i,j}\int_{B_{\tau}\setminus B_{\rho}}u_iu_j\partial_i(W_R)_j\dx\\
&\lesssim \left(\int_{B_\tau\setminus B_\rho}|u|^3\dx\right)^{2/3}\left(\int_{B_\tau\setminus B_R}|\nabla W_R|^3\dx\right)^{1/3}\\
& \lesssim (\tau-\rho)^{-1}\left(\int_{B_\tau\setminus B_\rho}|u|^3\dx\right)^{2/3}\left (\int_{B_\tau\setminus B_R}|u|^3\dx\right)^{1/3}\\
& \lesssim (\tau-\rho)^{-1}\int_{B_\tau\setminus B_\rho}|u|^3\dx,
\end{align*}
and similarly, we get
\begin{align*}
\int_{\bbR^3}(b \cdot \nabla)b\cdot W_R\dx 
& \lesssim (\tau-\rho)^{-1}\left(\int_{B_\tau\setminus B_\rho}|b|^3\dx\right)^{2/3}\left (\int_{B_\tau\setminus B_R}|u|^3\dx\right)^{1/3}\\
&\lesssim (\tau-\rho)^{-1}\int_{B_\tau\setminus B_\rho}|b|^3\dx+ (\tau-\rho)^{-1}\int_{B_\tau\setminus B_\rho}|u|^3\dx,
\end{align*}
where we have used Young's inequality. Furthermore, note that
\[
\int_{\bbR^3}\nabla \times w\cdot W_R\leq \|
w\|_{L^2(B_\tau)}\|\nabla W_R\|_{L^2(B_\tau)}\leq
\frac{1}{4}\|w\|^2_{L^2(B_\tau)}+C(\tau-\rho)^{-2}\|u\|^2_{L^2(B_\tau\setminus
 B_\rho)}.
\]
Finally, by H\"older's inequality and Young's inequality, we note that
\[\int_{\R^3}(u\cdot b)(b\cdot\nabla)\varphi^2_R\dx\lesssim \int_{B_{\tau}\setminus B_{\rho}}|u||b|^2|\nabla\varphi|\dx\lesssim(\tau-\rho)^{-1}\int_{B_{\tau}\setminus B_{\rho}}\left(|u|^3+|b|^3\right)\dx.\]
Altogether, we obtain from \eqref{main_ineq} that
\begin{align*}
&\hspace{4mm}\int_{B_{\rho}}\left(|\nabla u|^2+|\nabla b|^2+|\nabla w|^2\right)\dx+\int_{B_{\rho}}|w|^2\dx \\
& \lesssim (\tau-\rho)^{-2}\int_{B_\tau\setminus B_\rho}\left(|u|^2+|b|^2+|w|^2\right)\dx+(\tau-\rho)^{-1}\int_{B_\tau\setminus B_\rho}\left(|u|^3+|b|^3+|w|^3\right)\dx\\
&\hspace{4mm}+\frac{1}{4}\int_{B_\tau}\left(|\nabla u|^2+|\nabla b|^2+|\nabla w|^2\right)\dx+\frac{1}{4}\int_{B_\tau}|w|^2\dx.
\end{align*}
Then by Lemma \ref{key_2}, we conclude that
\[\hspace{-58mm}\int_{B_\rho}\left(|\nabla u|^2+|\nabla b|^2+|\nabla w|^2\right)\dx+\int_{B_{\rho}}|w|^2\dx\]
\vspace{-2mm}
\begin{equation}\label{main_ineq_2}
\lesssim (\tau-\rho)^{-2}\int_{B_\tau\setminus B_\rho}\left(|u|^2+|b|^2+|w|^2\right)\dx+(\tau-\rho)^{-1}\int_{B_\tau\setminus B_\rho}\left(|u|^3+|b|^3+|w|^3\right)\dx.
\end{equation}

Before proceeding further, let us briefly describe the strategy of the proof. We set $\tau=2\rho$ for convenience and we shall first show

\begin{equation}\label{main-est_1}
    \rho^{-1} \int _{B_{2\rho}\setminus B_\rho} (| u |^3 + | b |^3+ | w |^3) \,\mathrm{d}x \to 0 \qquad \mbox{as} \,\,\, \rho\rightarrow \infty.
\end{equation}
For the remaining part in \eqref{main_ineq_2}, by the H\"{o}lder's
inequality, we note that
\begin{equation*}
    \rho^{-2} \int_{B_{2\rho}\setminus B_\rho} (|u| ^2 + |b| ^2+ |w| ^2) \dx \lesssim \rho^{-\frac 13}
    \left( \rho^{-1} \int_{B_{2\rho}\setminus B_\rho} (|u|^3 + |b|^3+ |w|^3) \dx \right)^{\frac 23}.
\end{equation*}
Hence, if we show
\begin{equation}\label{main-est_2}
\left( \rho^{-1} \int_{B_{2\rho}\setminus B_\rho} (|u|^3 + |b|^3+ |w|^3) \dx \right)<C
\end{equation} 
for some constant $C>0$, we can get
\begin{equation}\label{main-est_3}
    \rho^{-2} \int_{B_{2\rho}\setminus B_\rho} (|u| ^2 + |b| ^2+ |w| ^2) \dx \to 0 \qquad \mbox{as} \,\,\, \rho \to \infty.
\end{equation}
Due to \eqref{main-est_1} and \eqref{main-est_3}, we have from \eqref{main_ineq_2} that
\begin{equation*}
\int_{B_\rho} (|\nabla u|^2 + |\nabla b|^2+ |\nabla w|^2) \dx \to 0 \qquad \mbox{as} \qquad \rho \rightarrow \infty,
\end{equation*}
which implies that $u$, $b$ and $w$ must be constants. Thanks to
\eqref{main-est_1}, we finally conclude that $u \equiv b \equiv w \equiv0$.

As described above, we first aim to prove \eqref{main-est_1}. Recall that $R>\rho>R/4>1$ and set $\psi= \xi_{\rho,R} - \xi_{\rho/4, R/4}$. By Lemma~\ref{LEM3}, we have from Young's inequality that
\begin{align*}
\int_{B_R} |\psi^3 u |^3 \dx 
& \lesssim R \| \psi \nabla u \|_{L^2(B_R)}^{\frac {18}{s+6}} + R^{4-\frac 3s} (R - \rho)^{-3}\\
&\hspace{4mm}+R \left( (R-\rho)^{-2} R^{\frac {11}6 - \frac 1s} \| \psi \nabla u \|_{L^2(B_R)} + R^{\frac {11}3 - \frac 2s} (R- \rho)^{-4} \right)^{\frac
    9{s+6}}\\
    &\lesssim R \| \psi \nabla u \|_{L^2(B_R)}^{\frac {18}{s+6}} +  R^{4-\frac 3s} (R - \rho)^{-3}+ R \big( \| \psi \nabla u \|_{L^2(B_R)}^2 + \rho^{\frac
{11}3 - \frac 2s} (R- \rho)^{-4} \big)^{\frac 9{s+6}}.
\end{align*}
By taking $\rho = 2r$ and $R = 4r$ for $r>1$, we deduce that
\begin{equation}\label{step_mid_1}
r^{-1} \int_{B_{2r} \setminus B_r} | u |^3 \dx \lesssim \|
\nabla u \|_{L^2(B_{4r}\setminus B_{r/2})}^{\frac {18}{s+6}} + r^{-\frac 3s}.
\end{equation}
Similarly, we can also obtain
\begin{equation}\label{step_mid_2}
r^{-1} \int_{B_{2r} \setminus B_r} \left(| b |^3+| w |^3\right) \dx
\lesssim \left(\| \nabla b \|_{L^2(B_{4r}\setminus B_{r/2})}^{\frac {18}{s+6}}+\| \nabla w
\|_{L^2(B_{4r}\setminus B_{r/2})}^{\frac {18}{s+6}}\right) + r^{-\frac 3s}.
\end{equation}

Next, we shall show that 
\begin{equation}\label{claim-second}
    \int_{\R^3} \left(|\nabla u|^2 + |\nabla b|^2+ |\nabla w|^2\right) \dx \leq C
\end{equation}
for some constant $C>0$. We set $R>\rho>1$ and $\overline{R}=(R+\rho)/2$. If we take $\varphi=\xi_{\rho,\overline{R}}$ as a cut-off function and proceed with the same argument used to derive \eqref{main_ineq_2}, we obtain that
\begin{align*}
& \hspace{4mm} \int_{B_\rho}\left(|\nabla u|^2+|\nabla b|^2+|\nabla w|^2\right)\dx\\
& \lesssim (R-\rho)^{-2}\int_{B_{\overline{R}}\setminus B_\rho}\left(|u|^2+|b|^2+|w|^2\right)\dx+(R-\rho)^{-1}\int_{B_{\overline{R}}\setminus B_\rho}\left(|u|^3+|b|^3+|w|^3\right)\dx.
\end{align*}
And then we set $\psi = \xi_{\overline{R}, R}$. Then as $\psi = 1$ on $B(\overline{R})$, we note that
\[
\int_{B_\rho} (|\nabla u|^2 + |\nabla b|^2+ |\nabla w|^2)\dx \leq C({\rm{I}}_1+{\rm{I}}_2),
\]
where
\[
{\rm{I}}_1 := (R-\rho)^{-2} \int_{B_R} \left(|\psi^2 u|^2 + |\psi^2 b|^2+
|\psi^2 w|^2\right)\dx,
\]
and
\[
{\rm{I}}_2 := (R-\rho)^{-1} \int_{B_R} \left(|\psi^3 u|^3 + |\psi^3 b|^3+
|\psi^3 w|^3\right) \dx.
\]
Then by Lemma \ref{LEM3} together with the assumption \eqref{AS}, and Young's inequality, we have that
\begin{align*}
(R-\rho)^{-1} \int_{B_R} |\psi^3 u|^3 \dx 
& \leq CR (R-\rho)^{-1} \| \psi \nabla u \|_{L^2(B_R)}^{\frac {18}{s+6}} + C R^{4-\frac 3s} (R - \rho)^{-4}\\
& \hspace{4mm} + C R(R-\rho)^{-1} \big( (R-\rho)^{-1} \| \psi^2 u \|_{L^2(B_R)} \big)^{\frac{18}{s+6}}\\
&\leq \frac{1}{4} \| \psi \nabla u \|_{L^2(B(R))}^2 + C R^{\frac {s+6}{s-3}}(R-\rho)^{-\frac {s+6}{s-3}} + {\rm{I}}_1.
\end{align*}

In the same way, if we proceed with the above argument for $b$ and $w$, we obtain
\begin{equation*}
{\rm{I}}_2 \leq \frac{1}{4} \int_{B_R} \left(|\nabla u|^2 + |\nabla b|^2+ |\nabla w|^2\right) \dx +C R^{\frac {s+6}{s-3}}(R-\rho)^{-\frac {s+6}{s-3}} + {\rm{I}}_1.
\end{equation*}
Next, for ${\rm{I}}_1$, using \eqref{C1}, we have by Young's inequality that
\begin{align*}
    (R-\rho)^{-2} \int_{B_R} |\psi^2 u|^2 \dx &\leq CR^{\frac {11}6 - \frac 1s} (R-\rho)^{-2} \| \psi \nabla u \|_{L^2(B_R)}
    + CR^{\frac {11}3 - \frac 2s} (R- \rho)^{-4} \\
    &\leq C R^2 (R-\rho)^{-2} \| \psi \nabla u \|_{L^2(B_R)} + R^4 (R-
    \rho)^{-4}\\
    &\leq\frac{1}{4} \| \psi \nabla u \|_{L^2(B_R)}^2 + C R^4 (R- \rho)^{-4}.
\end{align*}
If we use the same method with $b$ and $w$, it follows that
\[
{\rm{I}}_1 \leq \frac{1}{4} \int_{B_R} \left(|\nabla u|^2 + |\nabla b|^2+ |\nabla w|^2\right)\dx + C R^4 (R- \rho)^{-4}.
\]
Collecting the estimates for ${\rm{I}}_1$ and ${\rm{I}}_2$ yields
\[
\int_{B_\rho} \left(|\nabla u|^2 + |\nabla b|^2+ |\nabla w|^2\right) \dx \leq \frac{1}{2} \int_{B_R} \left(|\nabla u|^2 + |\nabla b|^2 + |\nabla w|^2\right) \dx + C R^{\frac {s+6}{s-3}}(R-\rho)^{-\frac {s+6}{s-3}},
\]
where we have used the facts $R(R-\rho)^{-1}>1$ and $3<s\leq 6$.
Then applying Lemma \ref{key_2} gives us the estimate
\[
\int_{B_\rho} \left(|\nabla u|^2 + |\nabla b|^2+ |\nabla w|^2\right) \dx \leq C R^{\frac {s+6}{s-3}}(R-\rho)^{-\frac {s+6}{s-3}}.
\]
If we take $R = 2\rho$ and let $\rho \rightarrow \infty$, we obtain \eqref{claim-second}, and hence from \eqref{step_mid_1} and \eqref{step_mid_2}, we conclude that \eqref{main-est_1} holds.

It remains to show \eqref{main-est_2}.

By direct computation, we observe that
\begin{align*}
r^{-1}\int_{B_r}\left(|u|^3+|b|^3+|w|^3\right)\dx
&=\sum_{j=1}^{\infty}2^{-j}(2^{-j}r)^{-1}\int_{B_{2^{-(j-1)}r}\setminus B_{2^{-j}r}}\left(|u|^3+|b|^3+|w|^3\right)\dx\\
& \leq \sup_{1/2 \leq \rho \leq r/2}\rho^{-1} \int _{B_{2\rho}\setminus B_\rho} \left(| u |^3 +| b |^3+| w |^3\right) \dx\\
&\hspace{4mm} + \int_{B_1} \left(| u |^3 +| b |^3+| w |^3\right) \dx.
\end{align*}
Therefore from \eqref{main-est_1}, we have \eqref{main-est_2}, and consequently, we deduce that the convergence \eqref{main-est_3} holds.

Now we are ready to conclude $u \equiv b \equiv w \equiv 0$. From \eqref{main_ineq_2} together with \eqref{main-est_1} and \eqref{main-est_3}, we have
\[
\int_{B_\rho} \left(|\nabla u|^2 + |\nabla b|^2+ |\nabla w|^2\right) \dx \to 0 \qquad \mbox{as} \,\,\, \rho \rightarrow \infty,
\]
which means that $u$, $b$ and $w$ are constants. By \eqref{main-est_1}, we finally obtain that $u \equiv b \equiv w \equiv0$.

\section{Proof of Theorem~\ref{thm2}}
From \eqref{main_ineq} in the proof
in Theorem \ref{thm1}, we know that
\[\hspace{-80mm}\int_{B_\rho}\left(|\nabla u|^2+|\nabla b|^2+|\nabla w|^2\right)\dx\]
\vspace{-2mm}
\begin{equation}\label{thm2_1}
\lesssim (\tau-\rho)^{-2}\int_{B_\tau\setminus B_\rho}\left(|u|^2+|b|^2+|w|^2\right)\dx+(\tau-\rho)^{-1}\int_{B_\tau\setminus B_\rho}\left(|u|^3+|u||b|^2+|u||w|^2\right)\dx
\end{equation}
In order to deal with the right-hand side of \eqref{thm2_1}, let us consider two cases as follows:

{\textbf{(Case 1)}} $p \in [3,\frac{9}{2})$: Then we have by H\"older's inequality,
\begin{align*}
\int_{B_\rho}\left(|\nabla u|^2+|\nabla b|^2+|\nabla w|^2\right)\dx 
&\lesssim(\tau-\rho)^{-2}\|(u, b, w)\|^2_{L^2(B_{2\rho})} + (\tau-\rho)^{-1}\|(u, b, w)\|^3_{L^3(B_{2\rho})}\\
& \lesssim \rho^{1-\frac{6}{p}}\|(u, b, w)\|^2_{L^p(B_{2\rho})} + \rho^{2-\frac{9}{p}}\|(u, b, w)\|^3_{L^p(B_{2\rho})},
\end{align*}
where we have chosen $\tau=2\rho$. Hence if we let $\rho\rightarrow\infty$, we have
\begin{equation}\label{mid_con_1}
\int_{\R^3}\left(|\nabla u|^2+|\nabla b|^2+|\nabla w|^2\right)\dx =0.
\end{equation}

{\textbf{(Case 2)}} $p \in [2,\frac{18}{5})$: For this case, let us consider a non-negative cut-off function $\theta(x) = \xi_{\tau, 2R}(x)-\xi_{\frac{R}{2},\rho}(x)$ with $\|\nabla \theta\|_{L^\infty}\lesssim \max\{\frac{1}{\rho-\frac{R}{2}}, \frac{1}{2R-\tau}\}$. Note that by the interpolation inequality,
\begin{align*}
\|w\|^2_{L^4(B_\tau\setminus B_\rho)}\lesssim\|w\theta\|^2_{L^4(\R^3)}
& \lesssim \|w\theta\|_{L^2(\R^3)}^{\frac{1}{2}}\|\nabla(w\theta)\|_{L^2(\R^3)}^{\frac{3}{2}}\\
& \lesssim \|w\theta\|^{\frac{1}{2}}_{L^2(\R^3)}\left(\|(\nabla w)\theta\|^{\frac{3}{2}}_{L^2(\R^3)}+\|w(\nabla\theta)\|^{\frac{3}{2}}_{L^2(\R^3)}\right)\\
& \lesssim \|w\|^{\frac{1}{2}}_{L^2(B_{2R}\setminus B_{R/2})}\left(\|\nabla w\|^{\frac{3}{2}}_{L^2(B_{2R}\setminus B_{R/2})}+\frac{1}{R^{\frac{3}{2}}}\|w\|^{\frac{3}{2}}_{L^2(B_{2R}\setminus B_{R/2})}\right)\\
& \lesssim \|w\|^{\frac{1}{2}}_{L^2(B_{2R}\setminus B_{R/2})}\|\nabla w\|^{\frac{3}{2}}_{L^2(B_{2R}\setminus B_{R/2})}+\frac{1}{R^{\frac{3}{2}}}\|w\|^{2}_{L^2(B_{2R}\setminus B_{R/2})}.\\
\end{align*}
And thus, by H\"older's inequality and Young's inequality,
\begin{align*}
(\tau-\rho)^{-1} \int_{B_\tau\setminus B_\rho} |u||w|^2\dx
& \lesssim (\tau-\rho)^{-1}\|u\|_{L^2(B_\tau\setminus B_\rho)}\|w\|^2_{L^4(B_\tau\setminus B_\rho)}\\
& \lesssim(\tau-\rho)^{-1}\|u\|_{L^2(B_{2R})}\left(\|w\|^{\frac{1}{2}}_{L^2(B_{2R})}\|\nabla w\|^{\frac{3}{2}}_{L^2(B_{2R})}+\frac{1}{R^{3/2}}\|w\|^{2}_{L^2(B_{2R})}\right)\\
& \lesssim(\tau-\rho)^{-4}\|(u,w)\|^{6}_{L^2(B_{2R})}+\frac{1}{4}\|\nabla
w\|^{2}_{L^2(B_{2R})}+\frac{(\tau-\rho)^{-1}}{R^{3/2}}\|(u,w)\|^{3}_{L^2(B_{2R})}.
\end{align*}
Therefore, we have
\begin{align*}
\int_{B_\rho}\left(|\nabla u|^2+|\nabla b|^2+|\nabla w|^2\right)\dx
& \leq \frac{1}{4} \int_{B_{2R}}\left(|\nabla u|^2+|\nabla b|^2+|\nabla w|^2\right)\dx+\frac{C}{(\tau-\rho)^2}\|(u,b,w)\|^2_{L^2(B_{2R})}\\
& \hspace{4mm} +\frac{C}{(\tau-\rho)^4}\|(u,b,w)\|^6_{L^2(B_{2R})} + \frac{C}{(\tau-\rho)R^{\frac{3}{2}}}\|(u,b,w)\|^3_{L^2(B_{2R})}\\
& \leq \frac{1}{4} \int_{B_{2R}}\left(|\nabla u|^2+|\nabla b|^2+|\nabla w|^2\right)\dx+CR^{1-\frac{6}{p}}\|(u,b,w)\|^2_{L^p(B_{2R})}\\
& \hspace{4mm} +CR^{5-\frac{18}{p}}\|(u,b,w)\|^6_{L^p(B_{2R})} + CR^{2-\frac{9}{p}}\|(u,b,w)\|^3_{L^p(B_{2R})}.
\end{align*}

If we apply Lemma \ref{key_2} with $f(r):=\int_{B_r} (|\nabla u|^2 + |\nabla b|^2+ |\nabla
w|^2) \,\mathrm{d}x$, we get
\begin{align*}
\int_{B_{R/2}}\left(|\nabla u|^2+|\nabla b|^2+|\nabla w|^2\right)\dx 
&\lesssim R^{1-\frac{6}{p}}\|(u,b,w)\|^2_{L^p(B_{2R})}+R^{5-\frac{18}{p}}\|(u,b,w)\|^6_{L^p(B_{2R})}\\
&\hspace{4mm}  + R^{2-\frac{9}{p}}\|(u,b,w)\|^3_{L^p(B_{2R})}.
\end{align*}
If we let $R \rightarrow \infty$ we can immediately find that
\[
\int_{\R^3} \left(|\nabla u|^2 + |\nabla
b|^2+ |\nabla w|^2\right) \dx=0,
\]
which implies that $u$, $b$ and $w$ are a constant vectors in
$R^3$. Since $u,b,w \in L^p(\R^3)$ for $p \in [2,\frac{9}{2})$, we conclude that $u \equiv b \equiv w \equiv 0$ in $\R^3$.

\section{Proof of Theorem~\ref{thm3}}

We first recall the fact that any continuous functions vanishing at infinity
must be bounded; thus, from the assumption of the theorem, we have $u,b,w \in L^\infty(\R^3)$. From
\eqref{main_ineq_2} in the proof in Theorem \ref{thm1}, we know that

\[\hspace{-80mm}\int_{B_\rho}\left(|\nabla u|^2+|\nabla b|^2+|\nabla w|^2\right)\dx\]
\vspace{-2mm}
\begin{equation}\label{thm3_1}
\lesssim \frac{1}{(\tau-\rho)^2}\int_{B_\tau\setminus B_\rho}\left(|u|^2+|b|^2+|w|^2\right)\dx+\frac{1}{(\tau-\rho)}\int_{B_\tau\setminus B_\rho}\left(|u|^3+|b|^3+|w|^3\right)\dx.
\end{equation}
In order to control the right-hand side in \eqref{thm3_1}, let us
consider it in two cases as follows:

{\textbf{(Case 1)}} $p\in [1,3/2)$: We note that with the choice $\tau = 2R$ and $\rho = R$,
\[
(\tau-\rho)^{-2} \int_{B_\tau\setminus B_\rho} |u| ^2\dx\lesssim (\tau-\rho)^{-2} \|u\|_{L^\infty(\R^3)}\int_{B_\tau\setminus B_\rho} |u| \dx \lesssim R^{1-\frac{3}{p}}\|u\|_{L^p(B_{2R})}.
\]
And thus, with the same arguments for $b$ and $w$, we get
\[
(\tau-\rho)^{-2}\int_{B_\tau\setminus B_\rho}\left(|u|^2+|b|^2+|w|^2\right)\dx \lesssim R^{1-\frac{3}{p}}\|(u,b,w)\|_{L^p(B_{2R})}.
\]

For the second term on the right-hand side of \eqref{thm3_1},
\[
(\tau-\rho)^{-1} \int_{B_\tau\setminus B_\rho} |u| ^3\dx\lesssim (\tau-\rho)^{-1} \|u\|^2_{L^\infty(\R^3)}\int_{B_\tau\setminus B_\rho} |u| \dx \lesssim R^{2-\frac{3}{p}}\|u\|_{L^p(B_{2R})},
\]
and similarly for $b$ and $w$, we have
\begin{equation}\label{final_1}
(\tau-\rho)^{-1}\int_{B_\tau\setminus B_\rho}\left(|u|^3+|b|^3+|w|^3\right)\dx \lesssim R^{2-\frac{3}{p}}\|(u,b,w)\|_{L^p(B_{2R})}.
\end{equation}

{\textbf{(Case 2)}} $p\in [\frac{3}{2},\frac{9}{4})$: In this case, we note that
\[
(\tau-\rho)^{-2} \int_{B_\tau\setminus B_\rho} |u| ^2\dx\lesssim(\tau-\rho)^{-2} \|u\|^{1/2}_{L^\infty(\R^3)}\int_{B_\tau\setminus B_\rho} |u|^{\frac{3}{2}} \dx \lesssim R^{1-\frac{9}{2p}}\|u\|^{\frac{3}{2}}_{L^p(B_{2R})},
\]
and in the same way for $b$ and $w$, we have
\[
(\tau-\rho)^{-2}\int_{B_\tau\setminus B_\rho}\left(|u|^2+|b|^2+|w|^2\right)\dx \lesssim R^{1-\frac{9}{2p}}\|(u,b,w)\|^{\frac{3}{2}}_{L^p(B_{2R})}.
\]

For the second term, in a similar way, we also get,
\[
(\tau-\rho)^{-1} \int_{B_\tau\setminus B_\rho} |u| ^3\dx\lesssim (\tau-\rho)^{-1} \|u\|^{\frac{3}{2}}_{L^\infty(\R^3)}\int_{B_\tau\setminus B_\rho} |u|^{\frac{3}{2}} \dx \lesssim R^{2-\frac{9}{2p}}\|u\|^{\frac{3}{2}}_{L^p(B_{2R})},
\]
which directly implies that
\begin{equation}\label{final_2}
(\tau-\rho)^{-1}\int_{B_\tau\setminus B_\rho}\left(|u|^3+|b|^3+|w|^3\right)\dx \lesssim R^{2-\frac{9}{2p}}\|(u,b,w)\|^{\frac{3}{2}}_{L^p(B_{2R})}.
\end{equation}

Collecting \eqref{final_1} and \eqref{final_2}, we have
\[
\int_{B_{R}} \left(|\nabla u|^2 + |\nabla b|^2+ |\nabla w|^2\right)
\dx\lesssim\begin{cases} (R^{1-\frac{3}{p}} + R^{2-\frac{3}{p}} )\|(u,b,w)\|_{L^p(B_{2R})},\quad\, p \in [1, \frac{3}{2}),\\
(R^{1-\frac{9}{2p}} + R^{2-\frac{9}{2p}}
)\|(u,b,w)\|^{3/2}_{L^p(B_{2R})},\quad p \in [\frac{3}{2},
\frac{9}{4}).
\end{cases}
\]

If we let $R\rightarrow +\infty$ in the above estimate, we can immediately find that
\[
\int_{\R^3} \left(|\nabla u|^2 + |\nabla b|^2+ |\nabla w|^2\right)\,dx=0,
\]
which implies that $u$, $b$ and $w$ are a constant vectors in $\R^3$. As we know $u,b,w \in L^p(\R^3)$ for $p \in [1,\frac{9}{4})$, we finally conclude that $u \equiv b \equiv w \equiv 0$ in $\bbR^3$.

\section*{Acknowledgments}
Seungchan Ko was supported by INHA UNIVERSITY Research Grant and National Research Foundation of
Korea Grant funded by the Korean Government (RS-2023-00212227).
Jae-Myoung Kim was supported by National Research Foundation of
Korea Grant funded by the Korean Government
(NRF-2020R1C1C1A01006521).

\section*{Data Availability Statement}
The data that supports the findings of this study are available within the article.

\bibliography{references}
\bibliographystyle{abbrv}

\end{document}